\documentclass[11pt,a4paper]{amsart}
\usepackage{amsmath,amsfonts,amssymb,amsthm}
\usepackage{graphicx}
\usepackage{wrapfig}
\newtheorem {lemma}{Lemma}
\newtheorem {corollary}{Corollary}
\newtheorem {proposition}{Proposition}
\newtheorem {theorem}{Theorem}
\theoremstyle{definition}
\newtheorem {example}{Example}
\newtheorem {definition}{Definition}
\newtheorem {remark}{Remark}
\def\ed{\ensuremath{\leftrightarrow}}

\DeclareSymbolFont{rsfscript}{OMS}{rsfs}{m}{n}
\DeclareSymbolFontAlphabet{\mathrsfs}{rsfscript}

\DeclareMathAlphabet{\mathbbold}{U}{bbold}{m}{n}

\DeclareFontFamily{OMS}{rsfs}{\skewchar\font'177}
\DeclareFontShape{OMS}{rsfs}{m}{n}{%
      <5> rsfs5
      <6> <7> rsfs7
      <8> <9> <10> rsfs10
      <10.95> <12> <14.4> <17.28> <20.74> <24.88> rsfs10
      }{}

\def\calO{\mathrsfs{O}}

\newcommand{\dl}{\ensuremath{<_{dl}}}

\renewcommand{\L}{{\mathrsfs{L}\!ie}}
\renewcommand{\P}{\mathrsfs{P}}
\newcommand{\G}{\mathrsfs{G}\!erst}
\newcommand{\Com}{{\mathrsfs{C}\!om}}
\newcommand{\LL}{{\mathrsfs{L}\!ie_2}}

\newcommand{\PP}{{\mathrsfs{P}_2}}
\newcommand{\GG}{\mathrsfs{G}\!erst_2}
\newcommand{\A}{{\ensuremath{OS^+}}}
\renewcommand{\AA}{{\ensuremath{OS^+_2}}}
\newcommand{\B}{{\ensuremath{OS}}}
\newcommand{\BB}{{\ensuremath{OS_2}}}
\DeclareMathOperator{\ord}{ord}

\makeatletter
\let\@newpf\proof \let\proof\relax 
\newenvironment{proof}{\@newpf[\proofname]}{\qed\endtrivlist}
\def\hm#1{#1\nobreak\discretionary{}{\hbox{\m@th$#1$}}{}}
\makeatother

\newlength{\wrfwidth}

\def\wrapfig#1#2#3#4{%
\settowidth{\wrfwidth}{\includegraphics[scale=#3]{#2}}
\begin{wrapfigure}{#1}{\wrfwidth}
\vskip#4
\centering
\includegraphics[scale=#3]{#2}%
\end{wrapfigure}}

\title[Quadratic algebras related to the bi-Hamiltonian operad]{Quadratic algebras\\ related to the bi-Hamiltonian operad}
\author{Mikhail Bershtein}
\address{Independent University of Moscow, Bolshoj Vlasievsky per., 11, Moscow, Russia, 119002}
\email{mbersht@mail.ru}
\author{Vladimir Dotsenko}
\address{Independent University of Moscow, Bolshoj Vlasievsky per., 11, Moscow, Russia, 119002}
\email{dotsenko@mccme.ru}
\author{Anton Khoroshkin}
\address{Institute for Theoretical and Experimental Physics, Bolshaya Cheremush\-kinskaya, 25, Moscow, Russia, 117218}
\email{khorosh@itep.ru}
\date{}

\begin{document}
\begin{abstract}
We prove the conjectures on dimensions and characters of some quadratic algebras stated by B.L.Feigin. It turns out that these algebras are naturally isomorphic to the duals of the components of the bi-Hamiltonian operad.
\end{abstract}
\maketitle

\section{Introduction.}

All algebras and operads in this paper are defined over the field of rational numbers $\mathbb{Q}$.

\subsection{Summary of results.}
The following series of quadratic algebras are well known.

\begin{definition}
The \emph{Orlik--Solomon algebra} $\B(n)$ is a super-commutative associative algebra with odd generators
$x_{ij}$, $1\le i\ne j\le n$, and relations
\begin{gather*}
x_{ij}-x_{ji}=0, \\
x_{ij}^2=0, \\
x_{ij}x_{jk}+x_{jk}x_{ki}+x_{ki}x_{ij}=0.
\end{gather*}
The \emph{even Orlik--Solomon algebra} $\A(n)$ is a super-commutative associative algebra with even generators
$x_{ij}$, $1\le i\ne j\le n$, and relations
\begin{gather*}
x_{ij}+x_{ji}=0, \\
x_{ij}^2=0, \\
x_{ij}x_{jk}+x_{jk}x_{ki}+x_{ki}x_{ij}=0.
\end{gather*}
\end{definition}

These algebras are closely related to the type $A$ hyperplane arrangements. The Orlik--Solomon algebra $\B(n)$ is isomorphic to the cohomology algebra of the complement of the arrangement $A_{n-1}$ over the complex numbers (see \cite{Ar}, \cite{OS}), while the even Orlik--Solomon algebra $\A(n)$ gives a graded version of the algebra of locally constant functions on the complement of the arrangement $A_{n-1}$ over the real numbers (see \cite{GV}). 

\begin{theorem}[Arnold~\cite{Ar}, Mathieu~\cite{Mat}]\label{OS}
 $$
\dim\B(n)=\dim\A(n)=n!.
 $$
\end{theorem}

Another interpretation of the Orlik--Solomon algebras was found in \cite{Coh}. One of the central results in \cite{Coh} is the computation of the homology for the (topological) operad of little discs. This homology (being an operad itself) turns out to be isomorphic to the Gerstenhaber operad (which is a graded version of the Poisson operad). Since components of the operad of little discs are homotopically equivalent to the complements of complex hyperplane arrangements, cohomology algebras of these components are isomorphic to Orlik--Solomon algebras. Thus there exists a cooperadic structure on the collection of Orlik--Solomon algebras. The even Orlik--Solomon algebras are also equipped with a cooperadic structure (and the dual operad is the Poisson operad). We show how these structures themselves can be used to compute the dimensions and monomial bases for these algebras.

The main results of our paper deal with the following generalisations of these algebras suggested by B.L.Feigin. These algebras are related to some other quadratic algebras introduced by A.N.Kirillov \cite{ANK}. Both our algebras and the ``diagonal coinvariants'' of Haiman \cite{Ha} coincide with certain quotients of Kirillov's algebras, and this could possibly lead to some new interpretation of diagonal coinvariants.

\begin{definition}
The \emph{double Orlik--Solomon algebra} $\BB(n)$ is an associative super-commutative algebra with
odd generators $x_{ij}$, $y_{ij}$, $1\le i\ne j\le n$, and relations
\begin{gather*}
x_{ij}-x_{ji}=y_{ij}-y_{ji}=0,\\
x_{ij}x_{jk}+x_{jk}x_{ki}+x_{ki}x_{ij}=0,\\
x_{ij}y_{jk}+x_{jk}y_{ki}+x_{ki}y_{ij}+y_{ij}x_{jk}+y_{jk}x_{ki}+y_{ki}x_{ij}=0,\\
y_{ij}y_{jk}+y_{jk}y_{ki}+y_{ki}y_{ij}=0,\\
x_{ij}^2=x_{ij}y_{ij}=y_{ij}^2=0.
\end{gather*}
The \emph{double even Orlik--Solomon algebra} $\AA(n)$ is an associative super-commutative algebra with
even generators $x_{ij}$, $y_{ij}$, $1\le i\ne j\le n$, and relations
\begin{gather*}
x_{ij}+x_{ji}=y_{ij}+y_{ji}=0,\\
x_{ij}x_{jk}+x_{jk}x_{ki}+x_{ki}x_{ij}=0,\\
x_{ij}y_{jk}+x_{jk}y_{ki}+x_{ki}y_{ij}+y_{ij}x_{jk}+y_{jk}x_{ki}+y_{ki}x_{ij}=0,\\
y_{ij}y_{jk}+y_{jk}y_{ki}+y_{ki}y_{ij}=0,\\
x_{ij}^2=x_{ij}y_{ij}=y_{ij}^2=0.
\end{gather*}
\end{definition}

The main result of this paper is

\begin{theorem}\label{OS2}
 $$
\dim\BB(n)=\dim\AA(n)=(n+1)^{n-1}.
 $$
\end{theorem}

\begin{remark}
There exists a neater way to describe the relations for $\AA(n)$. Let $z_{ij}=\lambda x_{ij}+\mu y_{ij}$. Then the relations mean that for each choice of coefficients $\lambda$ and $\mu$ the elements $z_{ij}$ satisfy the relations of~$\A(n)$. The analogous description for $\BB(n)$ fails: due to skew-commutativity, the relations $x_{ij}y_{ij}=0$ do not follow. Remarkably, the algebra obtained from the algebra $\B(n)$ in a similar way seems to be very interesting as well. For example, its dimension is conjectured to be $2^n(n+1)^{n-2}$ (A.N.Kirillov), which is also the conjectured dimension of the space of ``three-diagonal harmonics'' \cite{Ha}.
\end{remark}

Unfortunately, we do not know of any relationship between these algebras and hyperplane arrangements. Nevertheless, collections of these algebras form cooperads, and we can use operadic formalism to prove our theorem. It turns out that there exists a natural pairing between even Orlik--Solomon algebras and components of the bi-Hamiltonian operad. This operad (introduced in \cite{DK}) corresponds to algebraic structures arising on the ring of functions on a manifold equipped with two compatible Poisson brackets. Thus we also obtain additional structures on the bi-Hamiltonian operad. Our results provide another interpretation of the Hopf operad structure on this operad analogous to a well known Hopf operad structure on the Poisson operad. (Existence of a Hopf structure means that the tensor product of two bi-Hamiltonian algebras naturally possesses a structure of a bi-Hamiltonian algebra.) 

\begin{remark}
All known proofs for ordinary Orlik--Solomon algebras do not work in our case. The problem is that our algebras are not of PBW type (\cite{PP}): in fact, they are not Koszul; the algebra $\AA(3)$ has Betti number $b_{4,5}=2$,  and the algebra $\BB(3)$ has Betti number $b_{3,4}=1$ (both were computed by A.\,N.\,Kiril\-lov \cite{ANK} using the \texttt{bergman} computer algebra system). Thus there is no straightforward reduction of the problem to a combinatorial problem on monomial algebras. 
\end{remark}

\subsection{Plan of the paper.}The paper is organised as follows. 
In Section \ref{Operad} we recall the necessary definitions. The details can be found in \cite{DK}.
In Section \ref{Graph} we fix some notation for the various types of graphs used in our proofs.
To make the proof more transparent, we include a new proof of Theorem \ref{OS} (where we make use of known facts about the Poisson operad) in Section 4.
In Section \ref{Bihamiltonian} we prove Theorem \ref {OS2}. In fact, in both cases we prove only the part concerning the even Orlik--Solomon algebras; the other case is similar, and proofs are omitted.
First of all, we prove upper bounds on the dimensions of our algebras. To do this, we take a basis of the bi-Hamiltonian operad and introduce a bijection between this basis and a certain set of monomials in our algebras, and then show that these monomials form a spanning set for our algebras. (These monomials are indexed by some special kind of graphs.) We then prove a lower bound. We define a pairing between our algebras and components of the bi-Hamiltonian operad. It remains to prove that this pairing is non-degenerate. To do this, we introduce a linear ordering of our basis in the bi-Hamiltonian operad and the constructed spanning set of our algebras and then show that the corresponding pairing matrix is non-degenerate. [The proof for the double Orlik--Solomon algebras is essentially the same, but one has to replace the bi-Hamiltonian operad by the bi-Gerstenhaber operad.]

\subsection{Acknowledgements.}We are grateful to B.\,Feigin and A.\,N.\,Kirillov for numerous useful discussions. The second author is thankful to A.\,Szenes and M.\,Vergne for their interest in his work. We would like to thank the referees for several helpful remarks. We are also thankful to Alexander Frolkin for improving the language of our original text and for several typesetting hints. The work of the second author is partially supported by INTAS grant No.03-3350 and a LIEGRITS fellowship (contract No.MRTN-CT-2003-505078).  The work of the third author is partially supported by INTAS grant No.03-3350 and RFBR grant 04-02-16538.

\section{Definitions and notation.}\label{Operad}

We consider operads as certain functors from the category $Set_b$ of finite sets (with bijections as morphisms) to the category of vector spaces. Throughout the paper we take the space of multilinear elements in the free $\calO$-algebra whose generators are indexed by a finite set $I$ as a realisation of the space $\calO(I)$. For the ``standard'' finite set $[n]$, we call the vector space $\calO([n])$ the $n$th component of the operad $\calO$. In a similar way, we consider our algebras as functors from $Set_b$ to the category of algebras, allowing indices in the definition to range over a finite set $I$. (Thus it may be better to write $\A([n])$ instead of~$\A(n)$, but we prefer to keep the simpler notation.)

Let us remind the reader of the definitions of several operads. 

\subsection{The operads $\Com$, $\L$ and $\P$.}
We recall several standard definitions. The operad $\L$ is generated by a skew-symmetric binary operation $\{\cdot,\cdot\}$ with one quadratic relation: we wish the Jacobi identity~\eqref{Jacobi} to be satisfied in each algebra over this operad. Thus an algebra over this operad is a Lie algebra.
The operad $\Com$ is generated by a symmetric binary operation $\star$ with one quadratic relation, the associativity law for this operation. An algebra over this operad is an associative commutative algebra.

The Poisson operad $\P$ is generated by a symmetric operation $\star$ and a skew-symmetric operation $\{\cdot,\cdot\}$; the symmetric operation generates a suboperad of~$\P$ isomorphic to $\Com$, the skew-symmetric operation generates a suboperad isomorphic
to $\L$ and the relations between these operations mean that the skew-symmetric operation is a derivation of 
the symmetric operation (``the Leibniz rule for differentiating a product''):
 $$
\{a,b\star c\}=\{a,b\}\star c+b\star\{a,c\}.
 $$
An algebra over this operad is called a Poisson algebra; an example of such an algebra is the algebra of functions on a manifold with the ordinary product and with the bracket defined by a Poisson bivector field.

\subsection{The operads $\LL$ and $\PP$.}
The operad $\LL$ (also called the operad of two compatible brackets) is generated by two skew-symmetric operations (brackets) $\{\cdot,\cdot\}_1$ and $\{\cdot,\cdot\}_2$. The relations in this operad mean that any linear combination of these brackets satisfies the Jacobi identity.
It is equivalent to the following identities in each algebra over this operad: the Jacobi identity
\begin{gather}\label{Jacobi}
\{a,\{b,c\}_1\}_1+\{b,\{c,a\}_1\}_1+\{c,\{a,b\}_1\}_1=0,\\
\{a,\{b,c\}_2\}_2+\{b,\{c,a\}_2\}_2+\{c,\{a,b\}_2\}_2=0
\end{gather}
for each of the brackets and the six-term relation
\begin{multline}\label{6term}
\{a,\{b,c\}_1\}_2+\{b,\{c,a\}_1\}_2+\{c,\{a,b\}_1\}_2+\\
+\{a,\{b,c\}_2\}_1+\{b,\{c,a\}_2\}_1+\{c,\{a,b\}_2\}_1=0
\end{multline}
between the brackets.

The bi-Hamiltonian operad $\PP$ is generated by three operations, namely, two skew-symmetric operations
($\{\cdot,\cdot\}_1$ and $\{\cdot,\cdot\}_2$) and a symmetric operation ($\star$). The skew-commutative operations are two compatible Lie brackets (that is, relations \eqref{Jacobi}--\eqref{6term} hold), the commutative operation is an associative product and each of the brackets is a derivation of this product,
\begin{gather}
\{a,b\star c\}_1=\{a,b\}_1\star c+b\star\{a,c\}_1,\label{Leibniz}\\
\{a,b\star c\}_2=\{a,b\}_2\star c+b\star\{a,c\}_2.
\end{gather}

\subsection{The operads $\G$ and $\GG$.}

The Gerstenhaber operad $\G$ is a graded version of the Poisson operad; it is generated by a symmetric operation (product) $\cup$ of degree zero and a graded Lie bracket $[\cdot{,}\cdot]$ of degree~$-1$. The relations (for homogeneous elements of a graded algebra over this operad) are defined as follows.
\begin{itemize}
\item[(i)] super-commutativity of the product: $a\cup b=(-1)^{|a|\cdot |b|}b\cup a$,
\item[(ii)] super skew-commutativity of the bracket: $[a,b]=(-1)^{|a|\cdot |b|}[b,a]$,
\item[(iii)] associativity of the product: $a\cup(b\cup c)=(a\cup b)\cup c$,
\item[(iv)] graded Jacobi identity for the bracket: $$(-1)^{|a|(|c|-1)}[a,[b,c]]+(-1)^{|b|(|a|-1)}[b,[c,a]]+(-1)^{|c|(|b|-1)}[c,[a,b]]=0,$$
\item[(v)] graded Leibniz rule: $[a,b\cup c]=[a,b]\cup c+(-1)^{|b|\cdot |c|}[a,c]\cup b.$
\end{itemize}

The bi-Gerstenhaber operad $\GG$ is a graded version of the bi-Hamiltonian operad; it is generated by a symmetric operation (product) $\cup$ of degree zero and two graded Lie brackets $[\cdot{,}\cdot]_1$, $[\cdot{,}\cdot]_2$ of degree~$-1$.  The relations (for homogeneous elements of a graded algebra over this operad) are defined as follows.
\begin{itemize}
\item[(i)] super-commutativity of the product: $a\cup b=(-1)^{|a|\cdot |b|}b\cup a$,
\item[(ii)] super skew-commutativity of the brackets: $[a,b]_i=(-1)^{|a|\cdot |b|}[b,a]_i$, $i=1,2$,
\item[(iii)] associativity of the product: $a\cup(b\cup c)=(a\cup b)\cup c$,
\item[(iv)] graded Jacobi identity for the brackets: 
\begin{multline*}
(-1)^{|a|(|c|-1)}[a,[b,c]_i]_i+(-1)^{|b|(|a|-1)}[b,[c,a]_i]_i+\\
+(-1)^{|c|(|b|-1)}[c,[a,b]_i]_i=0, \quad i=1,2,
\end{multline*}
\item[(v)] graded compatibility condition for the brackets
\begin{multline*}
(-1)^{|a|(|c|-1)}[a,[b,c]_1]_2+(-1)^{|b|(|a|-1)}[b,[c,a]_1]_2+\\
+(-1)^{|c|(|b|-1)}[c,[a,b]_1]_2++(-1)^{|a|(|c|-1)}[a,[b,c]_2]_1+\\
+(-1)^{|b|(|a|-1)}[b,[c,a]_2]_1+(-1)^{|c|(|b|-1)}[c,[a,b]_2]_1=0,
\end{multline*}
\item[(vi)] graded Leibniz rule: $$[a,b\cup c]_i=[a,b]_i\cup c+(-1)^{|b|\cdot |c|}[a,c]_i\cup b, \quad i=1,2.$$
\end{itemize}

\begin{remark}
Due to the Leibniz rule, components of the Poisson (resp., bi-Hamiltonian) operad  are spanned by monomials that are products of Lie monomials (resp., monomials in two compatible brackets). The same holds for the Gerstenhaber and the bi-Gerstenhaber operads.
\end{remark}

\section{Combinatorial structures.}\label{Graph}

As the Lie operad is a quotient of a free operad with one generator, we can index the natural monomial spanning set for $\L(I)$ 
by planar binary trees 
whose leaves are 
labelled by $I$. In the case of the operad of two 
\wrapfig{r}{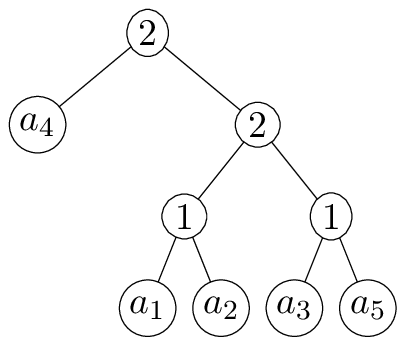}{1}{0cm}
compatible brackets, a spanning set is also indexed by binary trees, but these trees have an additional structure: the internal vertices are labelled by $\{1,2\}$
indicating which bracket is applied each time; see the picture for the tree corresponding to the monomial $\{a_4,\{\{a_1,a_2\}_1,\{a_3,a_5\}_1\}_2\}_2$. To any monomial which is a product of several monomials involving only brackets we assign the forest consisting of the corresponding binary trees.

Throughout the paper we use the combinatorics of graphs of two different types. The first type consists of (forests of) trees introduced above. They are of ``operadic'' nature and we call them $O$-trees when we want to avoid confusion. The second type of graph (``algebraic'' ones, referred to as $A$-graphs) is used to index spanning sets in our algebras. Vertices of these graphs are labelled by a finite set $I$. In the case of the algebras $\AA(I)$, the edges are 
labelled by the set $\{x,y\}$ indicating the generator. Let us describe precisely how monomials in our algebras correspond to graphs of this type.
\wrapfig{r}{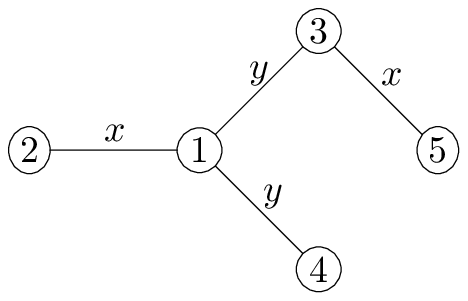}{1}{0.5cm}
\begin{definition} To a monomial from the algebra (of either type 
discussed above) with index set $I$, we assign a graph on $I$ in the following way. For each generator that appears in the monomial, let us connect $i$ and $j$ with an edge (and label this edge by $x$ or $y$ according to the type of the generator in the case of the algebras $\AA$; see the picture for the graph corresponding to the monomial $x_{12}y_{13}y_{14}x_{35}$). The converse is also clear: each graph on $I$ of this type gives rise to a monomial in our algebra. For a graph $\mathcal{G}$, we denote the corresponding monomial by $m_\mathcal{G}$.
\end{definition}

We will prove that for our algebras any monomial corresponding to a graph with a cycle vanishes, and so we only need to consider $A$-graphs which are forests of a certain type.

\section{The case of the Poisson operad.}\label{Poisson}

This section contains a new proof of previously known results on the Poisson operad. This proof provides a simple model of a more complicated proof for the bi-Hamiltonian operad and is organised similarly. 

\subsection{A monomial basis for the Poisson operad.} The results of this section are not new; they are included to make the proofs more structured and self-contained. 

We start with the operad $\L$. It is known that $\dim\L(n)=(n-1)!$ (for an operadic proof of this fact, see \cite{DK}). So to obtain a monomial basis in $\L(n)$, we need to find a spanning set of monomials of cardinality $(n-1)!$ in this space. The following lemma is classical.

\begin{lemma}\label{liebas}
The elements $\{a_{\sigma(1)},\{a_{\sigma(2)},\ldots,\{a_{\sigma(n-1)},a_n\}\ldots\}\}$ with $\sigma\in S_{n-1}$ span the space $\L(n)$.
\end{lemma}

\begin{proof}
The Jacobi identity, written in the form 
$$\{\{a,b\},c\}=\{a,\{b,c\}\}+\{b,\{a,c\}\},$$
allows us to decrease the degree (number of generators used) of the left argument of each bracket. Hence, by induction on the degree of the left argument, we can rewrite any monomial as a linear combination of monomials $\{a_i,m\}$. Using the skew symmetry of the bracket, we can assume at each step that the generator $a_n$ stays in the right argument, and the lemma follows by induction on $n$.
\end{proof}

Denote the set of all elements from Lemma~\ref{liebas} by $\mathfrak{B}_\L(n)$. It follows that this set is a basis for the component $\L(n)$. Since the Poisson operad is obtained from the operads $\Com$ and $\L$ by a distributive law, we have an isomorphism of~$\mathbb{S}$-modules $\P\simeq\Com\circ\L$ (see \cite{M}), and thus we have

\begin{theorem}
The set $\mathfrak{B}_\P(n)$ of monomials $c_1\star \ldots\star c_k$,
where $c_j\in \mathfrak{B}_\L(A_j)$ (for partitions $[n]=A_1\sqcup \ldots\sqcup A_k$ with $\max A_1<\ldots<\max A_k$), is a basis of~$\P(n)$.
\end{theorem}

\begin{corollary}\label{n-fact}
The dimension of the space $\P(n)$ is equal to $n!$.
\end{corollary}

\begin{proof}
Since the $\mathbb{S}$-modules $\P$ and $\Com\circ\L$ are isomorphic, the exponential generating series $f_\P(x)$ is equal to the composition 
$$f_\Com\circ f_\L=\exp(-\ln(1-x))-1=\frac{x}{1-x},$$ 
which is exactly the exponential generating function for $\{n!\}$. 

From the combinatorial point of view we can say that $(n-1)!$ is the number of cycles on $[n]$, while the composition with $\Com$ corresponds to a combinatorial structure obtained by partitioning a set and taking a cycle on each part of the partition. These structures are in one-to-one correspondence with permutations (decompositions into disjoint cycles), so the number of objects of this type is equal to $n!$.
\end{proof}

\subsection{An upper bound for $\dim\A(n)$.}

Here we prove that $n!$ is an upper bound for the dimension of~$\A(n)$. Moreover, we prove the following more general statement.

\begin{proposition}\label{upperbound}
Let $\A(n,\gamma)$ be an associative commutative algebra with generators
$x_{ij}$, $1\le i\hm\ne j\le n$, and relations
\begin{gather*}
x_{ij}+x_{ji}=0, \\
x_{ij}^2=0, \\
\gamma_{ijk} x_{ij}x_{jk}+\gamma_{jki} x_{jk}x_{ki}+\gamma_{kij} x_{ki}x_{ij}=0
\end{gather*}
for some $\gamma_{abc}\in\mathbb{Q}$. If $\gamma_{abc}\ne0$ for each $a,b,c\in [n]$, we have $\dim\A(n,\gamma)\le n!$.
\end{proposition}

\begin{proof}
To construct a spanning set for our algebra, consider the set $\mathfrak{F}_\A(n)$ of all monomials $m_\mathcal{G}$ corresponding to graphs $\mathcal{G}$ of the following type. Each connected component of~$\mathcal{G}$ (which is a graph on a subset $\pi\subset[n]$) is obtained from a linear reordering of~$\pi$ such that the maximal element of~$\pi$ remains maximal: from each element, edges go to its neighbours with respect to the ordering. 

To prove that these monomials span the algebra, we begin with the following lemma.

\begin{lemma}
Under the conditions of Proposition~\ref{upperbound}, the monomial $m_\mathcal{G}$ vanishes unless $\mathcal{G}$ is a forest.
\end{lemma}

\begin{proof}
We need to prove that $m_\mathcal{G}=0$ if $\mathcal{G}$ contains a cycle. 
We prove this by induction on the number of edges of the shortest cycle in $\mathcal{G}$.
If this number is equal to $2$, $\mathcal{G}$ has a multiple edge, and we use the relation $x_{ij}^2=0$. If the length is greater than $2$, consider two adjacent edges $i\leftrightarrow j$ and $j\leftrightarrow k$ of the cycle. Since $\gamma_{ijk}x_{ij}x_{jk}=-\gamma_{jki}x_{jk}x_{ki}-\gamma_{kij}x_{ki}x_{ij}$, the monomial $m_\mathcal{G}$ is equal to a linear combination of two monomials, each corresponding to a graph that has a cycle of smaller length. 
\end{proof}

Thus to prove the proposition, we need to rewrite monomials corresponding to a tree as linear combinations of the monomials introduced above. Consider a tree $T$. Let the vertex with maximal label be the root of~$T$. Consider the induced orientation of the tree (the edge $i\leftrightarrow j$ is oriented from $i$ to $j$ if the path from the root to $j$ goes via $i$). We prove our statement by induction on the distance $D$ from the root to the nearest vertex having at least two outgoing (directed) edges and, for a fixed $D$, on the number $N$ of outgoing edges of this vertex. Namely, we describe a procedure that rewrites a monomial as a linear combination of monomials corresponding to graphs with larger $D$ or the same $D$ and smaller $N$. It is clear that iterations of such a procedure lead to a linear combination of the monomials introduced above.

So fix $D\ge0$ and let $i$ be a vertex at distance $D$ having at least two outgoing edges $i\ed j$ and $i\ed k$. Using the relation $$\gamma_{kij}x_{ki}x_{ij}=-\gamma_{ijk}x_{ij}x_{jk}-\gamma_{jki}x_{jk}x_{ki},$$ we replace the corresponding monomial by a linear combination of two monomials where $i$ has fewer outgoing edges.

The  number of elements in the spanning set can be computed in several ways. To make the proof similar to the case of the bi-Hamiltonian operad, we provide a bijection~$\eta$ between this set and the basis of the component of the Poisson operad. 

Notice that our bases of components of the Lie operad can be described recursively: each element is a bracket $\{a_i,b\}$, where $i<n$ and $b$ belongs to a basis for a ``smaller'' component. This recursive definition immediately gives us a recursive definition for a bijection: take the tree corresponding to $b$, take its vertex $j$ corresponding to the uppermost generator in $b$ (i.e., such that $b=\{a_j,c\}$) and connect this vertex to $i$ with an edge. Combinatorially, if we draw the trees in the plane in a way that at each vertex the subtree containing $n$ is the rightmost one, we see that leaves of our tree are arranged in a chain according to their height ($n$ is the lowest), and we join the corresponding vertices with edges according to the order in this chain (see the picture). 
 $$
\includegraphics[scale=1.0]{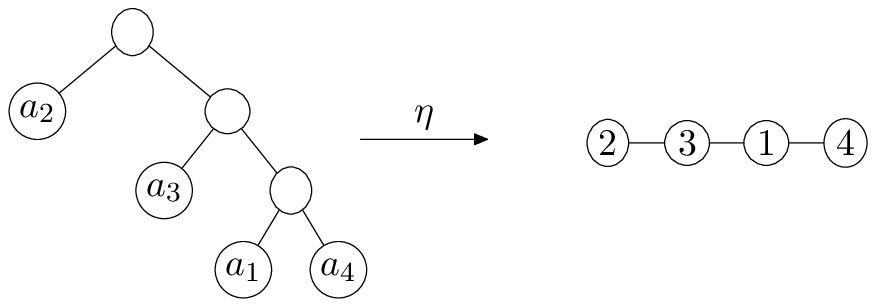}
 $$
Finally, we extend this correspondence on forests by taking the disjoint union of images of the trees in the forest.
\end{proof}

\subsection{A lower bound for $\dim\A(n)$.}

The upper bound $n!$ that we have obtained is not sharp for a generic algebra $\A(n,\gamma)$. To prove that it is sharp for $\A(n)$, we use operadic arguments. Namely, we introduce a pairing between the vector spaces $\A(n)$ and $\P(n)$ and then use this pairing to show that the spanning set from the previous section is in fact linearly independent. 

We define this pairing in two equivalent ways: one using a structure of a cooperad on the collection of algebras $\A(n)$ and the other one using a structure of a Hopf operad on $\P$. 

Let us begin with the cooperadic definition.

\begin{definition}
Take an element $\sharp\notin I$. Define an algebra homomorphism
$$\phi_{IJ}\colon \A(I\sqcup
J)\to \A(I\sqcup\{\sharp\})\otimes \A(J)$$
by the formulae
 $$
\phi_{IJ} (x_{ij})=
\left\{
\begin{aligned}
&x_{ij}\otimes 1, {\rm if\ }i,j\in I,\\
&1\otimes x_{ij}, {\rm if\ }i,j\in J,\\
&x_{i\sharp}\otimes 1, {\rm if\ }i\in I, j\in J.
\end{aligned}
\right.
 $$
\end{definition}

\begin{lemma}
For each $I$, $J$ the homomorphism $\phi_{IJ}$ is well defined (that is, the relations are mapped to zero).
\end{lemma}

The mappings $\phi_{IJ}$ obviously satisfy the coassociativity type relations for the cooperadic cocomposition, so the collection of the dual spaces $\{\A(n)^*\}$ acquires a structure of an operad.

\begin{lemma}
The binary operations of this operad satisfy the relations of~$\P$.
\end{lemma}

Thus there exists a mapping $\tau$ from $\P$ to the collection $\{\A(n)^*\}$. It maps the product $\star\in\P(2)$ to the linear function dual to $1\in\A(2)$ and maps the bracket $\{{,}\}$ to the linear function dual to $x_{12}\in\A(2)$. 

\begin{definition}
Define a pairing
 $$
\langle\cdot{,}\cdot\rangle\colon\P(n)\otimes \A(n)\to\mathbb{Q} 
 $$
by $\langle\alpha, m\rangle=(\tau(\alpha))(m)$.
\end{definition}

This definition can be explained in a combinatorial way. Take an $O$-forest $\mathcal{F}$ and an $A$-forest $\mathcal{G}$.
Define a partial mapping $f_{\mathcal{F},\mathcal{G}}$ from the set of edges of~$\mathcal{G}$ to the set of vertices of~$\mathcal{F}$ as follows. Given an edge $e=(i\leftrightarrow j)$, take the leaves $i$ and $j$ of~$\mathcal{F}$.
If they do not belong to the same component of~$\mathcal{F}$, then the mapping is not defined on~$e$. Otherwise, consider  for each of these leaves the path to the root of the corresponding connected component. Let $\kappa$ be the first common vertex of these paths. Then $f_{\mathcal{F},\mathcal{G}}(e)=\kappa$.

\begin{proposition}\label{pm1}
The value of the pairing $\langle\alpha_\mathcal{F},m_\mathcal{G}\rangle$ is nonzero iff the map $f_{\mathcal{F},\mathcal{G}}$ is a bijection between the set of all edges of~$\mathcal{G}$ and the set of all internal vertices of~$\mathcal{F}$. (In this case the value of the pairing is equal to $\pm1$.)
\end{proposition}

\begin{proof}
The mapping $\tau\colon\P(n)\to\A^*(n)$ can be thought of as follows. First of all, we apply the cocompositions in the cooperad dual to the Poisson operad to obtain an element of~$\P(2)^{\otimes (n-1)}$. This space is identified with the space $(\A(2)^*)^{\otimes (n-1)}$. Finally, this space is mapped to $\A(n)^*$ by the mapping which is dual to the corresponding cocomposition. Thus to each bracket in an operadic monomial the mapping $\tau$ assigns a generator in an algebraic monomial that corresponds to this bracket. One can easily see that this correspondence gives exactly the inverse to the mapping $f_{\mathcal{F},\mathcal{G}}$. The statement on the value of the pairing in the bijective case is clear. 
\end{proof}

Another way to define the pairing originates from the Hopf operad structure on $\P$ which we define now. Namely, we construct a coassociative morphism $\P\to\P\otimes\P$ (while all the properties of a Hopf operad hold, for us only the coassociativity property is important). For our purposes, the previous definition of the pairing is sufficient, but this Hopf approach to constructing the pairing can be used to obtain the exact formulae for the $\pm1$'s above.

\begin{definition}
Let $\Delta\colon\P(2)\to\P(2)\otimes\P(2)$
be defined by the formulae
\begin{gather*}
\Delta(\star)=\star\otimes\star,\\
\Delta(\{{,}\})=\{{,}\}\otimes\star+\star\otimes\{{,}\}.
\end{gather*}
Extend $\Delta$ by the rule $\Delta(a\circ b)=\Delta(a)\circ\Delta(b)$.
\end{definition}

The following lemma is straightforward.

\begin{lemma}
The (coassociative) morphism $\Delta$ is well-defined.
\end{lemma}

It follows that each vector space $\P(n)^*$ has the structure of an associative algebra.

Now we define the pairing. To make the definition compact, we use the following notation for the operations:
\begin{gather*}
\varepsilon\colon a_1,\ldots,a_n\mapsto a_1\star a_2\star\ldots\star a_n,\\
\alpha_{ij}\colon a_1,\ldots,a_n\mapsto
\{a_i, a_j\}\star a_1\star\ldots\star \hat{a}_i\star\ldots\star\hat{a}_j\star\ldots\star a_n.
\end{gather*}

\begin{definition}
Define the pairing 
 $$
\langle\cdot{,}\cdot\rangle\colon\P(n)\otimes \A(n)\to\mathbb{Q} 
 $$
as follows. Let $\langle\varepsilon, 1\rangle\hm=\langle\alpha_{ij}, x_{ij}\rangle\hm=1$, and let 
$\langle\alpha_{\mathcal{F}},m_{\mathcal{G}}\rangle=0$ for all other pairs $(\mathcal{F},\mathcal{G})$ where $\mathcal{G}$ has at most one edge. Extend this pairing to all $m_{\mathcal{G}}$, putting $\langle\lambda, gh\rangle=\langle\Delta(\lambda),g\otimes h\rangle$ (by
definition, the value of the pairing on tensor products is equal to the product of the values of the pairings on the factors).
\end{definition}

\begin{lemma}
This pairing is well-defined, i.e., it vanishes on all the
(operadic and algebraic) relations.
\end{lemma}

\begin{proposition}
$\dim\A(n)\ge n!$.
\end{proposition}

\begin{proof}
We begin with a definition of a linear ordering of the spanning set $\mathfrak{F}_\A(n)$ (and thus a linear ordering of~$\mathfrak{B}_\P(n)$ via the bijection $\eta$). After that we show that the matrix of the pairing between $\mathfrak{B}_\P(n)$ and $\mathfrak{F}_\A(n)$ is non-degenerate: for our orderings it turns out to be upper triangular (and in fact diagonal) with nonzero entries on the diagonal.

Now we  define a linear ordering of~$\mathfrak{F}_\A(n)$. Fix some ordering $<_p$ on the set of partitions $[n]=A_1\sqcup A_2\sqcup\ldots\sqcup A_k$ of~$[n]$.

Take two monomials $m_{\mathcal{G}_1}$ and $m_{\mathcal{G}_2}$ such that each $\mathcal{G}_i$ is a forest of trees that are chains having the maximal element among the leaves. For each of these forests, its connected components define a partition of~$[n]$. We say that the first element is greater than the second if the corresponding partition is greater with respect to the ordering $<_p$. Hence it remains to define the ordering for the forests such that the corresponding partitions are actually the same partition $[n]=A_1\sqcup\ldots\sqcup A_k$. Assume that the subsets $A_l$ are reordered in such a way that $\max A_1<\max A_2<\ldots<\max A_k$. 

We define the ordering for two trees on the same set of vertices, and then compare the trees on the subset $A_1$, the trees on the subset $A_2$,\ldots, the trees on the subset $A_k$. If at some of the comparison steps the corresponding trees are not equal, we take the first such step and say that the monomial having the greater tree on this step is greater than the other one.

So we take two trees $T$ and $T'$ on the same (ordered) vertex set $A$. Each of these trees is a chain of vertices starting at $\max A$, and we just compare the elements of these chains lexicographically: if for $T$ the vertex connected with $\max A$ is greater then for $T'$ then $T'<T$, if these vertices coincide, we compare their neighbours etc. 

It remains to prove that the matrix of the pairing between $\mathfrak{B}_\P(n)$ and $\mathfrak{F}_\A(n)$ is non-degenerate. Namely, we prove that $\langle\alpha_{\mathcal{F}},m_\mathcal{G}\rangle=0$ if $\mathcal{F}\hm<\eta^{-1}(\mathcal{G})$ (and $\langle\alpha_{\mathcal{F}},m_\mathcal{G}\rangle=\pm1$ if $\mathcal{F}=\eta^{-1}(\mathcal{G})$). Note that in this case we actually have $\langle\alpha_{\mathcal{F}},m_\mathcal{G}\rangle=0$ if $\mathcal{F}\ne\eta^{-1}(\mathcal{G})$, but we prefer the weaker statement to make the similarity with the proof for the bi-Hamiltonian operad more transparent.

First of all, from the combinatorial definition of our pairing it is evident that it vanishes if the partitions of~$[n]$ defined by the forests $\mathcal{F}$ and $\mathcal{G}$ are distinct. So we can concentrate on the case of coinciding partitions $[n]= A_1\sqcup\ldots\sqcup A_k$. As above, we can assume without loss of generality that $\max A_1<\max A_2<\ldots<\max A_k$, and it is enough to prove that $\langle\alpha_{\mathcal{T}'},m_T\rangle=0$ for an operadic tree $\mathcal{T'}$ and an algebraic tree $T$ (on the same vertex set) such that $\eta(\mathcal{T}')<T$. Let us denote $T'=\eta(\mathcal{T}')$. To make the notation simpler, we assume that the common vertex set of~$T$ and $T'$ is $[n]$.

Indeed, let us start computing the value of the pairing. We want the mapping $f_{\mathcal{T}',T}$ to be bijective. This means that the parent vertex $v$ of the leaf $n$ in $\mathcal{T}'$ corresponds to some edge of~$T$. The only possible edge is the edge connecting $n$ and its left sibling $i$. Now we take the parent vertex of~$v$ and use the same argument: if there exists an edge corresponding to this vertex then the left sibling of~$v$ is connected in $T$ with either $n$ or $i$~--- but by definition $n$ has only one neighbour, and we already know that this neighbour is $i$. So there is an edge $i\leftrightarrow j$ in $T$ etc.
\end{proof} 

The upper and the lower bound together imply the following result:
 
\begin{theorem}
\begin{enumerate}
\item For each $n$ the algebras $\A(n)$ and $\P(n)^*$ (with the product given by the Hopf structure on the Poisson operad) are isomorphic. 
\item The dimension of~$\A(n)$ is $n!$. 
\item The basis in this algebra is naturally indexed by forests from $\mathfrak{F}_\A(n)$.
\end{enumerate}
\end{theorem}

\section{The case of the bi-Hamiltonian operad.}\label{Bihamiltonian}

\subsection{A monomial basis for the bi-Hamiltonian operad.}

One of the central results of the paper \cite{DK} is the dimension formula for the operad $\LL$. It states that $\dim\LL(n)=n^{n-1}$. Hence to introduce a monomial basis for this operad it is enough to find a spanning set having the same cardinality. 

\begin{definition}
Given a finite ordered set $A\hm=\{a_1,a_2,\ldots,a_n\}$, $a_1<a_2<\ldots<a_n$, define a family of monomials $\mathfrak{B}(A)$ in the free algebra with two compatible brackets generated by $A$ recursively as follows.
\begin{itemize}
\item For $A=\{a_1\}$, let $\mathfrak{B}(A)=\{a_1\}$.
\item For $n>1$, a monomial $b$ belongs to $\mathfrak{B}(A)$ if and only if it satisfies one of the two conditions:
\begin{enumerate}
\item $b=\{a_i, b'\}_1$, where $i<n$ and $b'\in\mathfrak{B}(A\setminus\{a_i\})$;
\item $b=\{b_1,b_2\}_2$, where $b_1\in\mathfrak{B}(A_1)$, $b_2\in\mathfrak{B}(A_2)$
for some $A_1 \sqcup A_2\hm=A$, $a_n\in A_2$, and either $b_1=a_i$ for some~$i$ or 
$b_1$ is a bracket of the first type of two monomials belonging to bases for some subsets of~$A_1$.
\end{enumerate}
\end{itemize}
\end{definition}

\begin{lemma}
$\mathfrak{B}_\LL(A)$ spans $\LL(A)$ and $|\mathfrak{B}_\LL(A)|=|A|^{|A|-1}$. 
\end{lemma}

\begin{proof}
Let us begin with the spanning property. Here the argument is similar to the case of the Poisson operad. Namely, we use the Jacobi identity for the first bracket and the six-term relation to rewrite each monomial as a linear combination of monomials, each either of the type $\{a_i,m\}_1$ with $i<n$ or of the type $\{m_1,m_2\}_2$. For monomials of the first type we just apply induction on $n$. As for the second type of monomials, we again assume that on each step $a_n$ belongs to the right argument, and the only case which does not allow us to use the inductive hypothesis is $m_1=\{m_1',m_2'\}_2$. In this case we should use the Jacobi identity for the second bracket to decrease the degree (number of generators involved) of the left argument, and the spanning property follows.

Let $\beta_n=|\mathfrak{B}_\LL([n])|$. Moreover, for $i=1,2$ let $\beta_{i,n}=|\mathfrak{B}^i_\LL([n])|$, 
where $\mathfrak{B}^1_\LL([n])$ (respectively, $\mathfrak{B}^2_\LL([n])$) is the set of all monomials described under item (1) (respectively, (2)) of the above definition. (We set $\beta_{1,1}=1$, $\beta_{2,1}=0$.) Let us obtain the recurrence relations for these sequences.

The first condition implies that
 $
\beta_{1,n+1}=n\beta_n,
 $
while the second one gives the relation
 $
\beta_{2,n+1}=\sum_{k=1}^n\binom{n}{k}\beta_{1,k}\beta_{n+1-k}
 $
(we choose $k$ indices to be involved in $b_1$; notice that $a_{n+1}$ is forbidden, so the factor is exactly $\binom{n}{k}$).

Let us rewrite these formulae as
\begin{gather*}
\frac{(n+1)\beta_{1,n+1}}{(n+1)!}=\frac{n\beta_n}{n!},\\
\frac{(n+1)\beta_{2,n+1}}{(n+1)!}=
\sum_{k=1}^n\frac{\beta_{1,k}}{k!}\frac{(n+1-k)\beta_{n+1-k}}{(n+1-k)!}. 
\end{gather*}
Using the exponential generating functions 
 $$
\beta(t)\hm=\sum_{l\ge1}\frac{\beta_lt^l}{l!}, \quad \beta_{i}(t)\hm=\sum_{l\ge1}\frac{\beta_{i,l}t^l}{l!},
 $$
we can encode the previous formulae via differential equations:
\begin{gather}
\beta'_1(t)-1=t\beta'(t),\label{first}\\
\beta'_2(t)=\beta_1(t)\beta'(t).\label{second}
\end{gather}
Substituting $\beta_2(t)=\beta(t)-\beta_1(t)$ into \eqref{second}, we get 
$\beta'(t)-\beta'_1(t)=\beta_1(t)\beta'(t)$,
which is equivalent to
 $$
\beta'(t)=\frac{\beta'_1(t)}{1-\beta_1(t)}.
 $$
Integrating the latter equality, we have $\beta(t)=-\ln(1-\beta_1(t))$,
and so $\exp(\beta(t))=\frac{1}{1-\beta_1(t)}$.
Rewrite this formula as $\beta_1(t)=1-\exp(-\beta(t))$ and substitute it into \eqref{first}.
We have $\beta'(t)\exp(-\beta(t))-1=t\beta'(t)$, and so   
 $$
\beta'(t)=\exp(\beta(t))+t\beta'(t)\exp(\beta(t)).
 $$
Integrating again, we have 
 $$
\beta(t)=t\exp(\beta(t)).
 $$ 
It is well known that the only solution of this functional equation is the generating function of~$\{n^{n-1}\}$.
\end{proof}

It follows that the sets $\mathfrak{B}_\LL([n])$ provide monomial bases for the components of the operad $\LL$. Since the bi-Hamiltonian operad is obtained from operads $\Com$ and $\LL$ by a distributive law, we have an isomorphism of~$\mathbb{S}$-modules $\PP\simeq\Com\circ\LL$ (see \cite{DK}), and thus we have

\begin{theorem}
The set $\mathfrak{B}_\PP(n)$ of monomials $c_1\star\ldots\star c_k$,
where $c_j\in \mathfrak{B}_\LL(A_j)$ (for all partitions  $[n]=A_1\sqcup \ldots\sqcup A_k$ with $\max A_1\le\ldots\le\max A_k$), is a basis of~$\PP(n)$.
\end{theorem}

This leads to the following ``combinatorial'' proof of the dimension formula for the bi-Hamiltonian operad (which is different from the proof given in \cite{DK}).

\begin{corollary}\label{biham-n-fact}
The dimension of the space $\PP(n)$ is $(n+1)^{n-1}$.
\end{corollary}

\begin{proof}
It is well known that $n^{n-1}$ is equal to the number of rooted trees on $[n]$. Taking the composition with $\Com$ corresponds to a combinatorial structure obtained by partitioning a set and taking a rooted tree on each part of the partition. This data is in one-to-one correspondence with planted forests on $[n]$, or, equivalently, trees on $[n+1]$, so the number of objects of this type is $(n+1)^{n-1}$.
\end{proof}

\subsection{An upper bound for $\dim\AA(n)$.}\label{Upper}

Here we obtain an upper bound for the dimension of~$\AA(n)$. Again, this upper bound holds for a generic algebra of this type.

\begin{proposition}
Let $\AA(n,\gamma)$ be an associative commutative algebra with generators
$x_{ij}$, $y_{ij}$, $1\le i\ne j\le n$, and relations
\begin{gather*}
x_{ij}+x_{ji}=y_{ij}+y_{ji}=0,\\
x_{ij}^2=x_{ij}y_{ij}=y_{ij}^2=0, \\
\gamma_{xx,ijk} x_{ij}x_{jk}+\gamma_{xx,jki} x_{jk}x_{ki}+\gamma_{xx,kij} x_{ki}x_{ij}=0,\\
\gamma_{xy,ijk} x_{ij}y_{jk}+\gamma_{xy,jki} x_{jk}y_{ki}+\gamma_{xy,kij} x_{ki}y_{ij}+
\phantom{\gamma_{xy,ijk} x_{ij}y_{jk}+\gamma_{xy,jki} x_{jk}y_{ki}}\\
\phantom{x_{ij}y_{jk}+\gamma_{xy,jki} x_{jk}y_{ki}}+\gamma_{yx,ijk} y_{ij}x_{jk}+\gamma_{yx,jki} y_{jk}x_{ki}+\gamma_{yx,kij} y_{ki}x_{ij}=0,\\
\gamma_{yy,ijk} y_{ij}y_{jk}+\gamma_{yy,jki} y_{jk}y_{ki}+\gamma_{yy,kij} y_{ki}y_{ij}=0
\end{gather*}
for some $\gamma_{\alpha\beta,abc}\in\mathbb{Q}$. 
If $\gamma_{\alpha\beta,abc}\ne0$ for each $a,b,c\in [n]$, $\alpha,\beta\in\{x,y\}$, we have $\dim\AA(n,\gamma)\le (n+1)^{n-1}$.
\end{proposition}

\begin{proof}
Let us define a bijection between our basis of~$\PP(n)$ and a certain set of monomials in $\AA(n)$, describe this set combinatorially, and then prove the spanning property.

Let us define a mapping from $\mathfrak{B}_\PP([n])$ to the set of~$A$-forests on $[n]$. Actually, we define some mapping $\psi$ from $\mathfrak{B}_\LL([n])$ to the set of~$A$-trees and then extend it on forests in an obvious way. 
Take an element $b\in\mathfrak{B}_\LL([n])$. Let us explain how to construct a tree $\psi(b)$ on the vertex set $[n]$ with $\{x,y\}$-labelling of the edges. The definition is recursive. For $n=1$ we do not have any choice (both of the sets consist of one element). Consider the case $n>1$. There are two different cases.

\begin{enumerate}
\item
$b=\{a_i,\{b_1,\{\ldots,\{b_k,b_{k+1}\}_2\ldots\}_2\}_2\}_1$, where each $b_s$ is either $a_{i_s}$ or $\{a_{i_s},b_s'\}_1$ for some $i_s$ and some basis element $b_s'$. Then the tree $\psi(b)$ is obtained by adding edges $i_1\leftrightarrow i_2$, $i_2\leftrightarrow i_3$,\ldots, $i_k\leftrightarrow i_{k+1}$ labelled~$y$ and an edge $i\leftrightarrow i_{k+1}$ labelled~$x$ to the trees $\psi(b_1)$, $\psi(b_2)$, \ldots, $\psi(b_{k+1})$.
\item $b=\{b',b''\}_2$, where $b'$ is either $a_i$ or $\{a_i,c'\}_1$ for some $a_i$ and some basis element $c'$, and
$b''$ is either $a_j$ or $\{a_j,c''\}_1$ or $\{\{a_j,c''\}_1, d''\}_2$ for some $j$ and some basis elements $c''$, $d''$. Then
the tree $\psi(b)$ is obtained by adding an edge $i\leftrightarrow j$ labelled~$y$ to the trees $\psi(b')$ and $\psi(b'')$.
\end{enumerate}

\begin{example}\label{Ex1}
$\psi(\{a_4,\{\{a_1,a_2\}_1,\{a_3,a_5\}_1\}_2\}_2)$ is the graph corresponding to the monomial $x_{12}y_{13}y_{14}x_{35}$ (see Section 3 for the corresponding picture).
\end{example}

The injectivity of~$\psi$ is clear. Let us describe the image $\psi(\mathfrak{B}_\LL([n]))$. 

Take any $A$-tree $T$. Let $n$ be its root; thus we obtain an orientation of edges of~$T$ in the usual way. Consider all trees $T$ satisfying the following condition: 
\begin{quote}
($\ast$)\ for each vertex of~$T$ there exists at most one outgoing (oriented) edge labelled~$x$ and
at most one outgoing (oriented) edge labelled~$y$.
\end{quote}
For each vertex $v$ of~$T$ we call any sequence $c$ of vertices $v_0=v$, $v_1$,\ldots, $v_k$, where all the oriented edges $v_i\rightarrow v_{i+1}$ are labelled by $x$, an $x$-chain starting at $v$ (in this case, we say that the length of~$c$ is equal to $k$). An $x$-chain is called maximal if it can not be extended to a longer $x$-chain. If the above condition is satisfied then for each vertex $v$ of~$T$ there exists exactly one maximal $x$-chain starting at $v$. We denote this chain by $c_x(v,T)$ and its length by $l_x(v,T)$. Maximal $y$-chains are defined in a similar way and the corresponding notation should be clear.

Any tree $T$ satisfying ($\ast$) defines a partition of~$[n]$ as follows. Consider the maximal $x$-chain $j_0(T)=n$, $j_1(T)$,\ldots, $j_k(T)$ starting at $n$ (here $k=l_x(n,T)$). Let us delete all the edges used in this chain. Connected components of the resulting graph are trees $T_0$, $T_1$,\ldots, $T_k$ on some subsets of $[n]$. These subsets form a partition $[n]=\xi_0(T)\sqcup\xi_1(T)\sqcup\ldots\sqcup\xi_k(T)$ (here $\xi_s(T)$ stands for the part containing $j_s(T)$). Take any $p=0,1,\ldots,k$ and consider the maximal $y$-chain $i_{0,p}(T)=j_p$, $i_{1,p}(T)$,\ldots, $i_{k_p,p}(T)$ starting at $j_p$ in $T_p$ (here $k_p=l_y(j_p,T)$). Let us delete all the edges used in that chain as well. Connected components of the resulting graph are trees $T_{i,p}$ ($i=1,\ldots,k_p$) on some subsets of $\xi_p(T)$. These subsets form a partition $$\xi_p(T)=\pi_{0,p}(T)\sqcup\ldots\sqcup\pi_{k_p,p}(T).$$
\begin{definition}\label{TreesOS}
Let $\mathfrak{T}_\AA([n])$ be the set of all $A$-trees satisfiying ($\ast$) and the following (recursive) conditions:
\begin{itemize}
\item[(i)] for all $p=0,\ldots, k$ and $s=1,\ldots,k_p$ the end of the maximal $x$-chain beginning at $i_{s,p}(T)$ is the maximal element of $\pi_{s,p}(T)$;
\item[(ii)] each of the trees $T_{i,p}$ obtained after deleting the edges of the maximal chains (as above) belongs to the set $\mathfrak{T}_\AA(\pi_{i,p})$.
\end{itemize}
\end{definition}

\begin{example}
The following tree $T$ belongs to $\mathfrak{T}_\AA(17)$.
 $$
\includegraphics{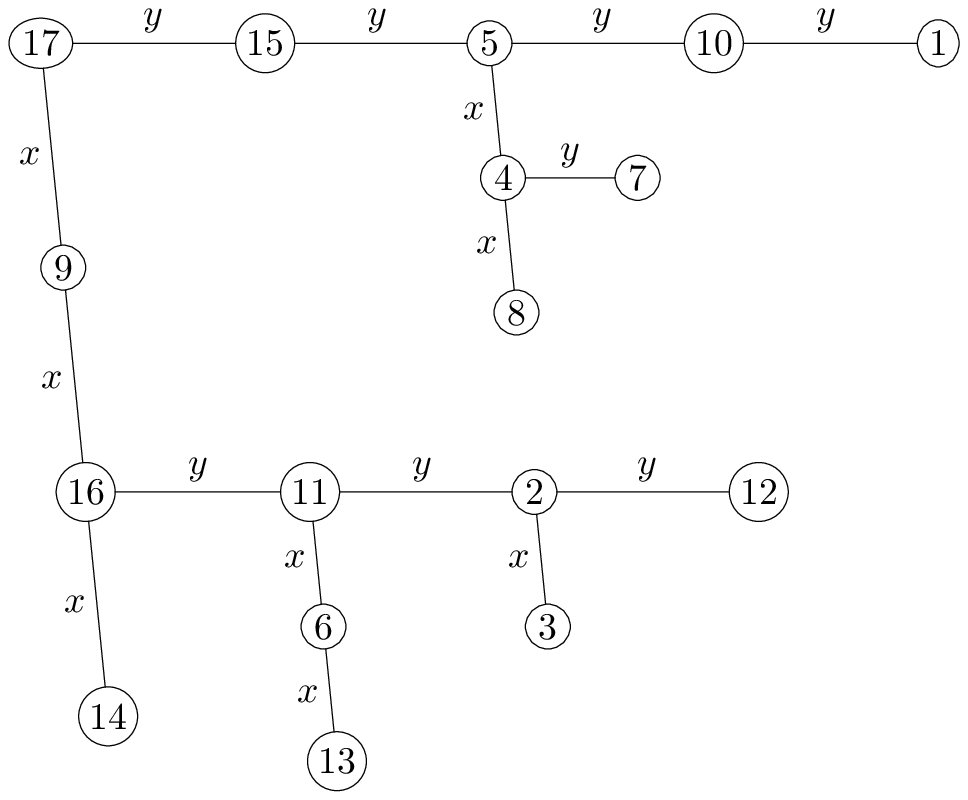}
 $$
In particular, 
 $$
l_x(17,T)=3, l_y(17,T)=4, c_x(5,T)=(5,4,8), \pi_{2,0}(T)=\{4,5,7,8\}.
 $$
\end{example}

\begin{definition}
Take a forest $\mathcal{F}$ on $[n]$ with $\{x,y\}$-labelled edges. Let $$[n]=\pi_1\sqcup\ldots\sqcup\pi_m$$ 
be the partition of $[n]$ defined by connected components $\mathcal{F}_1$, \ldots, $\mathcal{F}_m$ of $\mathcal{F}$. 
We denote by $\mathfrak{F}_\AA([n])$ the set of all forests $\mathcal{F}$ on $[n]$ such that for each $s$ the tree $\mathcal{F}_s$ belongs to $\mathfrak{T}_\AA(\pi_s)$.
\end{definition}

Let us explain why $\psi(\mathfrak{B}_\PP([n]))=\mathfrak{F}_\AA([n])$.
For $b\in\mathfrak{B}_\LL([n]))$, consider the $O$-tree corresponding to $b$. Take the unique path in this tree going from the leaf labelled~$n$ to the root. This path contains some internal vertices labelled~$1$ and some vertices labelled~$2$ which are split into connected components by vertices labelled~$1$. Notice that from the definition of the mapping $\psi$ it is clear that the vertices labelled~$1$ correspond to the vertices $j_0(T)=n$, \ldots, $j_m(T)$ of the maximal $x$-chain starting from $n$ in $\psi(b)$, and that the chain consisting of vertices labelled~$2$ between the $s$th and the $(s+1)$th vertex labelled $1$ corresponds to the maximal $y$-chain starting from $j_s(T)$. Now it is easy to see that the maximality conditions correspond to conditions $\max A_1<\max A_2$ from the definition of $\mathfrak{B}_\LL(A)$. We omit the further details.

We now want to prove that our monomials span $\AA(n)$.

\begin{lemma}
The monomials corresponding to forests span $\AA(n)$.
\end{lemma}

\begin{proof}
Let us prove that any monomial corresponding to a graph with cycles is equal to zero. We prove this by induction on the length of the shortest cycle. For a cycle of length two it follows literally from the defining relations. If a cycle has two adjacent edges labelled~$x$ (respectively, $y$) then we use the three-term $x$-relation (respectively, $y$-relation) to reduce the length of a cycle. The only remaining case is the case of edges with alternating labels. But then we use the six-term relation for two adjacent edges and obtain five monomials, four of which have a cycle of smaller
length and the fifth two adjacent edges with the same labels, and we are done.
\end{proof}

Now let us prove that for any forest $\mathcal{G}$ on $[n]$ the monomial $m_\mathcal{G}$ can be rewritten as a linear combination of monomials corresponding to forests from $\mathfrak{F}_\AA(n)$. We use induction on
$n$. Using the inductive hypothesis, we can immediately derive

\begin{lemma}
Each monomial $m_\mathcal{G}$ is equal to a linear combination of monomials $m_\mathcal{F}$ for forests $\mathcal{F}$ consisting of trees that satisfy the following (for some $k,l$):
\begin{itemize}
\item The vertex with the maximal label has $k$ outgoing edges labelled $x$ and $l$
outgoing edges labelled $y$.
\item Let us delete these outgoing edges. We will get trees $T^x_1$, \ldots, $T^x_k$, $T^y_1$, \ldots, $T^y_l$ on some subsets $\pi^x_1$, \ldots, $\pi^x_k$, $\pi^y_1$, \ldots, $\pi^y_l$. The condition is that each of the trees $T^y_j$ belongs to the corresponding set $\mathfrak{T}_\AA(\pi^y_j)$.
\end{itemize}
\end{lemma}

Let us prove that we can fulfil condition (i) of Definition \ref{TreesOS} for elements of~$y$-chains starting from $n$. It is enough to show this for a single chain. The proof is by decreasing induction on $k+l$ and for fixed $k+l$ by induction on $l$. In the first case, the induction base is $k+l=n-1$, and all the conditions are fulfilled. In the second case the induction base is $l=0$, so there are no $y$-chains and the conditions are fulfilled as well.

Consider one of the trees $T^y_p$. Let us take the maximal element $m^y_p$ of the set $\pi^y_p$. 
Denote by $v^y_p$ the endpoint of the $x$-chain starting at $m^y_p$. First of all, we prove that 
\begin{quote}
($\blacklozenge$) modulo the monomials corresponding to the trees satisfying the inductive hypothesis, our monomial is congruent to a monomial for which $n$ is connected to $v^y_p$ (and all other edges remain the same). 
\end{quote}
This monomial satisfies the maximality conditions for the corresponding $y$-chain starting from $n$. 

Consider the vertex $v$ of~$T^y_p$ that is connected to $n$. Let us take the path going from that vertex to $v^y_p$. Denote vertices of that path by $v_0=v^y_p$, \ldots, $v_r=v$. We prove~($\blacklozenge$) by induction on $r$ (base being obvious). Indeed, if the edge $v_r\ed v_{r-1}$ is labelled $y$, we use the relation $$\gamma_{yy,nv_rv_{r-1}}y_{nv_r}y_{v_rv_{r-1}}=
-\gamma_{yy,v_{r-1}v_rn}y_{nv_{r-1}}y_{v_{r-1}v_r}+\gamma_{yy,v_rnv_{r-1}}y_{nv_{r-1}}y_{nv_r},$$
where for the first monomial $r$ has decreased, while for the second $l$ (and so $k+l$) has increased. If the edge $v_r\ed v_{r-1}$ is labelled $x$, we use the relation
\begin{multline*}
\gamma_{yx,nv_rv_{r-1}}y_{nv_r}x_{v_rv_{r-1}}=-\gamma_{yx,v_{r-1}v_rn}y_{nv_{r-1}}x_{v_{r-1}v_r}
+\gamma_{yx,v_rnv_{r-1}}y_{nv_{r-1}}x_{nv_r}-\\-\gamma_{xy,nv_rv_{r-1}}x_{nv_r}y_{v_rv_{r-1}}-\gamma_{xy,v_{r-1}v_rn}x_{nv_{r-1}}y_{v_{r-1}v_r}+\gamma_{xy,v_rnv_{r-1}}x_{nv_{r-1}}y_{nv_r},
\end{multline*}
where for the first monomial $r$ has decreased, for the second and the fifth one $k+l$ has increased, for the third and the fourth monomials $k+l$ is the same and $k$ has increased, so $l$ has decreased (it may be that some of the monomials no longer belong to the corresponding set and we should rewrite them through the our spanning set, but that does not affect $k$ and $l$).

When we are done, we finish as follows. First of all, we straighten the ``skeleton'' consisting of all $y$-chains starting at $n$. This can be done in the same way as for $\A(n)$, and does not affect the maximality conditions involving elements of the chains. Then using the same ``straightening'' procedure, we replace each monomial by a linear combination of monomials corresponding to trees where $n$ has at most one outgoing $x$-edge (the subtree connected to $n$ by the $y$-edge remains the same). For each monomial of this kind we take the subtree $\widetilde{T}$ (on the vertex set $\widetilde{\pi}$) growing from that $x$-edge and rewrite the corresponding monomial as a linear combination of monomials from $\mathfrak{T}_\AA(\widetilde{\pi})$ (here the ordering of the set $\widetilde{\pi}$ differs from the induced ordering; we force the label of the vertex connected to $n$ to be the maximal element of this set). This can be done by the inductive hypothesis, and the proposition follows. 
\end{proof}

\subsection{A lower bound for $\dim\AA(n)$.}

The proof is completely analogous to the proof for the case of the algebras $\A(n)$.
We give similar definitions for the pairing and then define a linear ordering and check the triangularity property.

\begin{definition}
Take an element $\sharp\notin I$. Define an algebra homomorphism
$$\phi_{IJ}\colon \AA(I\sqcup J)\to \AA(I\sqcup\{\sharp\})\otimes \AA(J)$$
as follows.
 $$
\phi_{IJ} (x_{ij})=
\left\{
\begin{aligned}
&x_{ij}\otimes 1, {\rm if\ }i,j\in I,\\
&1\otimes x_{ij}, {\rm if\ }i,j\in J,\\
&x_{i\sharp}\otimes 1, {\rm if\ }i\in I, j\in J,
\end{aligned}
\right.
 $$
 $$
\phi_{IJ} (y_{ij})=
\left\{
\begin{aligned}
&y_{ij}\otimes 1, {\rm if\ }i,j\in I,\\
&1\otimes y_{ij}, {\rm if\ }i,j\in J,\\
&y_{i\sharp}\otimes 1, {\rm if\ }i\in I, j\in J.
\end{aligned}
\right.
 $$
\end{definition}

\begin{lemma}
For each $I$, $J$ the homomorphism $\phi_{IJ}$ is well-defined.
\end{lemma}

The mappings $\phi_{IJ}$ obviously satisfy the coassociativity type relations for the cooperadic cocomposition, so the collection of the dual spaces $\{\AA(n)^*\}$ acquires the structure of an operad.

\begin{lemma}
The binary operations of this operad satisfy the relations of~$\PP$.
\end{lemma}

Thus there exists a mapping $\tau$ from $\PP$ to the collection $\{\AA(n)^*\}$. It maps the product $\star\in\PP(2)$ to the linear function dual to $1\in\AA(2)$ and the brackets $\{{,}\}_1$, $\{{,}\}_2$ to the linear functions dual to $x_{12}, y_{12}\in\AA(2)$. 

\begin{definition}
Define a pairing
 $$
\langle\cdot{,}\cdot\rangle\colon\PP(n)\otimes \AA(n)\to\mathbb{Q} 
 $$
by $\langle\alpha, m\rangle=(\tau(\alpha))(m)$.
 \end{definition}

The combinatorial description of the pairing goes as follows. Take an $A$-forest $\mathcal{G}$ and an $O$-forest $\mathcal{F}$. Define a partial mapping $f_{\mathcal{F},\mathcal{G}}$ from the set of edges of~$\mathcal{G}$ to the set of vertices of~$\mathcal{F}$ as follows. Given an edge $e=(i\leftrightarrow j)$, take the leaves $i$ and $j$ of~$\mathcal{F}$. If they do not belong to the same component of~$\mathcal{F}$, then the mapping is not defined on~$e$. Otherwise, consider for each of these leaves the path to the root of the corresponding connected component. Let $\kappa$ be the first common vertex of these paths. If the label of~$\kappa$ is $1$ and the label of~$e$ is $y$, or, vice versa, the label of~$\kappa$ is $2$ and the label of~$e$ is $x$, the mapping is not defined on~$e$. Otherwise, $f_{\mathcal{F},\mathcal{G}}(e)=\kappa$. The following is analogous to Proposition~\ref{pm1}.

\begin{proposition}
The value of the pairing $\langle\alpha_{\mathcal{F}},m_\mathcal{G}\rangle$ is nonzero iff the map $f_{\mathcal{F},\mathcal{G}}$ is a bijection between the set of all edges of~$\mathcal{G}$ and the set of all internal vertices of~$\mathcal{F}$. (In this case the value of the pairing is equal to $\pm1$.)
\end{proposition}

Another way to define this pairing originates from the following Hopf operad structure on $\PP$. 

\begin{definition}
Let the morphism
$\Delta\colon\PP(2)\to\PP(2)\otimes\PP(2)$
be defined as follows:
\begin{gather*}
\Delta(\star)=\star\otimes\star,\\
\Delta(\{{,}\}_1)=\{{,}\}_1\otimes\star+\star\otimes\{{,}\}_1,\\
\Delta(\{{,}\}_2)=\{{,}\}_2\otimes\star+\star\otimes\{{,}\}_2.
\end{gather*}
Extend $\Delta$ by the rule $\Delta(a\circ b)=\Delta(a)\circ\Delta(b)$.
\end{definition}

\begin{lemma}
The (coassociative) morphism $\Delta$ is well-defined.
\end{lemma}

It follows that each vector space $\PP(n)^*$ possesses the structure of an associative algebra.

Fix the following notation for the operations:
\begin{gather*}
\varepsilon\colon a_1,\ldots,a_n\mapsto a_1\star a_2\star\ldots\star a_n,\\
\alpha_{ij}\colon a_1,\ldots,a_n\mapsto
\{a_i, a_j\}_1\star a_1\star\ldots\star \hat{a}_i\star\ldots\star\hat{a}_j\star\ldots\star a_n,\\
\beta_{ij}\colon a_1,\ldots,a_n\mapsto
\{a_i, a_j\}_2\star a_1\star\ldots\star \hat{a}_i\star\ldots\star\hat{a}_j\star\ldots\star a_n.
\end{gather*}

\begin{definition}
Define the pairing 
 $$
\langle\cdot{,}\cdot\rangle\colon\PP(n)\otimes \AA(n)\to\mathbb{Q} 
 $$
as follows. 
Let $\langle\varepsilon, 1\rangle\hm=\langle\alpha_{ij}, x_{ij}\rangle\hm=\langle\beta_{ij}, y_{ij}\rangle\hm=1$, and let 
$\langle\alpha_{\mathcal{F}},m_{\mathcal{G}}\rangle=0$ for all other pairs $(\mathcal{F},\mathcal{G})$ where $\mathcal{G}$ has at most one edge. Extend this pairing to all $m_{\mathcal{G}}$, putting $\langle\lambda, gh\rangle=\langle\Delta(\lambda),g\otimes h\rangle$.
\end{definition}

\begin{lemma}
This pairing is well defined, i.e., it vanishes on all the (operadic and algebraic) relations.
\end{lemma}

\begin{proposition}
We have $\dim\AA(n)\ge (n+1)^{n-1}$.
\end{proposition}

\begin{proof}
As in the Poisson case, we begin with the definition of a linear ordering of~$\mathfrak{F}_\AA(n)$ (and thus a linear ordering of~$\mathfrak{B}_\PP(n)$ via the bijection $\psi$). After that we show that the matrix of the pairing is upper triangular with nonzero entries on the diagonal and thus non-degenerate.

Take two monomials $m_{\mathcal{G}_1}$ and $m_{\mathcal{G}_2}$ such that each $\mathcal{G}_i$ is a forest of trees of~$\mathfrak{T}_\AA$-type. For each of these forests, its connected components define a partition of~$[n]$. We say that the first element is greater than the second if the corresponding partition is greater with respect to the ordering $<_p$. Hence it remains to define the ordering for the forests such that the corresponding partitions are actually the same partition $[n]=A_1\sqcup\ldots\sqcup A_k$. Assume that the subsets $A_l$ are reordered in such a way that $\max A_1<\max A_2<\ldots<\max A_k$. 

We define the ordering for two trees on the same set of vertices, and then compare the trees on the subset $A_1$, the trees on the subset $A_2$,\ldots, the trees on the subset $A_k$. If at some comparison steps the corresponding trees are
not equal, we take the first such step and say that the greater monomial is the one having the greater tree on this step.

The ordering is more complicated than in the case of the Poisson operad. For $n=2$ we just fix any of two possible orderings. Assume that $n>2$ and take two trees $T$ and $T'$ on the same vertex set $A$ of cardinality $n$. Now we introduce the recursive ordering procedure. 

Fix some ordering on the subsets extending the inclusion ordering and denote this ordering by $\prec$. Denote also by $\dl$ the degree-lexicographic ordering of finite sequences of integers ($s_1\dl s_2$ if either the length of~$s_1$ is less than the length of~$s_2$ or the lengths coincide, and $s_1$ is lexicographically less than $s_2$; for example, an empty chain is less than any other chain). 

In the following definition we use the notation introduced in Section \ref{Upper}.

\begin{itemize}
\item Let $T'<T$ if $c_x(n,T')\dl c_x(n,T).$ (Thus we can assume in further comparisons that $c_x(n,T)=c_x(n,T')$.)
\item Let $k=l_x(n,T)=l_x(n,T')$ and $(j_1,\ldots,j_k)=c_x(n,T)\hm=c_x(n,T')$. Let $T'<T$ if $\xi_k(T')\prec \xi_k(T).$ (Thus we can assume in further comparisons that $\xi_k(T)\hm=\xi_k(T')$.)
Let $\zeta_k(T)\subset[n]$ be the set of indices involved in $c_y(j_k,T)$.
\item Assume that the subtree $T_k$ contains no edges labelled $x$. (Thus $\zeta_k(T)=\xi_k(T)$.)
Let $T'<T$ if $$c_y(j_k,T')\dl c_y(j_k,T).$$ If the latter chains coincide, let $T'<T$ if $T'_-<T_-$,
where the trees $T_-$, $T'_-$ are obtained from $T$ and $T'$ by deleting the edge $j_{k-1}\leftrightarrow j_k$ and the subtree growing from $j_k$. (Thus if $T_k$ does not contain $x$-edges, the comparison procedure stops.)

If the subtree $T_k$ contains some edges labelled $x$, consider the maximal element $\mu(T)$ of~$\xi_k(T)\setminus\zeta_k(T)$. Suppose that it belongs to $\pi_{\alpha(T),k}$. Then $c_x(i_{\alpha(T),k},T)$ ends in $\mu(T)$ by definition.
\item Let $T'<T$ if $$c_x(i_{\alpha(T'),k},T')^{\mathrm{rev}}\dl c_x(i_{\alpha(T),k},T)^{rev},$$ where the superscript $\mathrm{rev}$ means the same sequence in reverse order
(the first element becomes the last etc.). (Thus we can assume in further comparisons that $c_x(i_{\alpha(T'),k},T')\hm=c_x(i_{\alpha(T),k},T)$.)
\item Let $T'<T$ if $\pi_{\alpha(T),k}\prec\pi_{\alpha(T'),k}.$ Notice that the the smaller subset corresponds to the greater tree. (Thus we can assume in further comparisons that $\pi_{\alpha(T),k}=\pi_{\alpha(T'),k}$.)
\item Let $T'<T$ if $T'_{\alpha(T'),k}<T_{\alpha(T),k}$, where $T_{\alpha(T),k}$ and $T'_{\alpha(T'),k}$ are the subtrees growing from $i_{\alpha(T),k}=i_{\alpha(T'),k}$ in $T$ and $T'$. If these subtrees coincide, let us delete them (leaving $i_{\alpha(T),k}=i_{\alpha(T'),k}$ in their place), obtaining trees $\overline{T}_{\alpha(T),k}$ and $\overline{T'}_{\alpha(T'),k}$.
\item Finally, let $T'<T$ if $\overline{T'}_{\alpha(T'),k}<\overline{T}_{\alpha(T),k}$.
\end{itemize}

It is clear that this comparison procedure leads to a certain linear ordering. Consider the matrix of the pairing (whose rows and columns are ordered according to our ordering). We  prove that $\langle\alpha_{\mathcal{F}},m_\mathcal{G}\rangle=0$ if $\mathcal{F}<\psi^{-1}(\mathcal{G})$. It is clear from the definition that $\langle\alpha_{\mathcal{F}},m_\mathcal{G}\rangle=\pm1$ if $\mathcal{F}=\psi^{-1}(\mathcal{G})$. It will follow that our matrix is triangular with $\pm1$'s on the diagonal and hence non-degenerate.

First of all, from the combinatorial definition of the pairing it is evident that the pairing vanishes if the partitions of~$[n]$ defined by the forests $\mathcal{F}$ and $\mathcal{G}$ are distinct. So we can concentrate on the case of coinciding partitions $[n]= A_1\sqcup\ldots\sqcup A_k$. As above, we can assume (without loss of generality) that $\max A_1<\max A_2<\ldots<\max A_k$, and it is enough to prove that $\langle\alpha_{\mathcal{T}'},m_T\rangle=0$ for an $O$-tree $\mathcal{T}'$ and an $A$-tree $T$ (on the same vertex set) such that $\psi(\mathcal{T}')<T$. Let us denote $T'=\psi(\mathcal{T}')$. Thus we have $T'<T$. To make the notation simpler, we assume that the common vertex set of~$T$ and $T'$ is $[n]$.

Let us trace the comparison procedure. Consider the path in $\mathcal{T}'$ from its leaf $n$ to the root. 
Let $\alpha_{\mathcal{T}'}=\{b_1,\{b_2,\dots\{b_r,a_n\}_{\epsilon_r}\ldots\}_{\epsilon_2}\}_{\epsilon_1}$, where
$b_{\alpha}\in\mathfrak{B}_\LL(\omega_{\alpha})$ and $\epsilon_\alpha\in\{1,2\}$. 
From the definition of~$\psi$ we see that 
 $$
r=l_x(n,T') + \sum_{p=0}^{l_x(n,T')} l_y(j_p(T'),T'), 
 $$
and the number of~$1$'s among $\epsilon_i$'s is exactly $l_x(n,T')$. Moreover, the partitions  $[n]=\bigsqcup\pi_{q,p}(T')$ and $[n]=\{n\}\sqcup\bigsqcup \omega_{\alpha}$ actually coincide (up to rearranging the parts). Let us denote by $V_n =V_n^1\sqcup V_n^2$ the subset of vertices (labelled by $1$, $2$ resp.) in $\mathcal{T}'$ that belong to the path from the root to $n$. 
Consider the numbering 
 $$
\ord\colon V_n\rightarrow\{1,\ldots,r\}
 $$ 
that corresponds to the order of vertices in the path from the root to $n$ 
(if $v$ is closer to the root than $w$, then $\ord(v)<\ord(w)$). For a vertex $v\in V_n^1$ with $\ord(v)=\alpha$, we have $|\omega_\alpha|=1$.

\begin{remark}\label{root_path}
Consider the map $f_{\mathcal{T}',T}$. This map takes the edge $i\ed j$ to $V_n$ if and only if $i\in \omega_{\alpha}$, $j\in \omega_{\beta}$ and $\alpha\neq\beta$. In this case, we have $\ord(f_{\mathcal{T}',T}(i\ed j))=\min(\alpha,\beta)$.
\end{remark}

Recall that if $\langle\alpha_{\mathcal{T}'},m_T\rangle\ne0$ then $f_{\mathcal{T}',T}$ should be a bijection. This simple observation is actually crucial in our proof.

The same arguments as for the Poisson operad show that all the images of edges of~$c_x(n,T')$ belong to $V_n^1$ and  $\ord(f_{\mathcal{T}',T}(j_p\leftrightarrow j_{p+1}))$ decreases as $p$ increases. Thus if $c_x(n,T')<_{dl}c_x(n,T)$ then $\langle\alpha_{\mathcal{T}'},m_T\rangle=0$.

Suppose the comparison procedure for $T$ and $T'$ stops on the second step. 
Since $f_{\mathcal{T}',T}$ should be a bijection, all images of~$T_k$'s edges should belong to the $O$-subtree growing from the vertex $v_k$ in $V_n$ with $\ord(v_k)=l_y(j_k(T'),T')+1$,  i.e. the vertex in $V_n^1$ that is closest to the root. This contradicts the condition $\xi_k(T')\prec\xi_k(T)$, and therefore $\langle\alpha_{\mathcal{T}'},m_T\rangle=0$.

The third and fourth steps are similar to the first one; we just repeat the arguments for the Poisson operad (in the case of the third step we apply the Poisson operad argument to $y$-chains starting at $j_k$ in both trees, while in the case of the fourth step we use that argument for $x$-chains starting at $\mu(T)$). 

Suppose that the comparison procedure stops on the fifth step. This means that there exists an element
$s\in \pi_{\alpha(T'),k}\setminus\pi_{\alpha(T),k}$. Suppose that $\mu(T)=\mu(T')$ belongs to the set of leaves of~$\omega_{\alpha_0}$. Consider the shortest path 
 $$
s_0:=\mu(T)\ed s_1 \ed\ldots\ed s_l = s
 $$
in the tree $T$ starting at $\mu(T)$ and ending at $s$. 
Clearly, its first segment $s_{0}\ed\ldots\ed s_{l_1}$ coincides with $c_x(i_{\alpha(T),k},T)^{rev}$.
Notice that the set of all edges $s_i\ed s_{i+1}$ that are taken to $V_n$ by $f_{\mathcal{T}',T}$ is not empty (for example, this is true for $i=l_1$). Let $l_1$, \ldots, $l_q$ be the numbers of these edges.  For each $p=1,\ldots,q-1$, images of all edges in the segment $s_{l_p+1}\ed\ldots\ed s_{l_{p+1}-1}$ belong to the $O$-subtree corresponding to the monomial $b_{\alpha_p}$ (for some vertex $\alpha_p\in V_n$). On the other hand, images of all edges of the segment 
$s_{l_q+1}\ed\ldots\ed s_l$ belong to the $O$-subtree $b_{\alpha_0}$.
Thus 
 $$
\ord(f_{\mathcal{T}',T}(s_{l_{p}}\ed s_{l_p+1}))=\min(\alpha_p,\alpha_{p+1})
 $$
(we set $\alpha_{q+1}:=\alpha_0$). Now the vanishing property for the fifth step of the ordering follows from the following elementary lemma.

\begin{lemma}
Let $d_0,d_1,\ldots,d_q\ (q>1)$ be any sequence of numbers with $d_0=d_q$ and 
let $e_k:=\min(d_{k-1},d_k)$, $k=1,\ldots,q$. Then the numbers $e_k$ cannot be all distinct.
\end{lemma}

Now we can use the inductive hypothesis in the last comparison step
to finish the proof of triangularity.
\end{proof}

The upper and the lower bound together imply our main result:
\begin{theorem}
\begin{enumerate}
\item For each $n$ the algebras $\AA(n)$ and $\PP(n)^*$ (with the product given by the Hopf structure on the bi-Hamiltonian operad) are isomorphic. 
\item The dimension of~$\AA(n)$ is equal to $(n+1)^{n-1}$.  
\item The basis in this algebra is naturally indexed by forests from $\mathfrak{F}_\AA(n)$.
\end{enumerate}
\end{theorem}

Using character formulae from \cite{DK}, one immediately has
\begin{corollary}
The $SL_2$-character of~$\AA(n)$ (for the $SL_2$-action arising from the natural action on the set of generators) is equal to
 $$
\prod_{k=1}^{n-1}(1+kq+(n-k)q^{-1}).
 $$
\end{corollary}
We refer the reader to \cite{DK} for the formulae for the $S_n\times SL_2$-characters. 

Using our monomial bases, one can easily compute the Hilbert series of our algebras. The proof is straightforward and therefore omitted.

\begin{proposition}
The Hilbert series of~$\AA(n)$ is equal to $(1+nt)^{n-1}$.
\end{proposition}

\end{document}